\documentclass[lettersize,journal]{IEEEtran}

\usepackage[justification=centering]{caption}
\usepackage{multirow}
\usepackage{pbox}
\usepackage{makecell}

\usepackage{amsmath,amsfonts}
\usepackage{array}
\usepackage{textcomp}
\usepackage{stfloats}
\usepackage{url}
\usepackage{verbatim}
\usepackage{graphicx}
\usepackage{subfigure}
\usepackage[justification=centering]{caption}
\usepackage{epstopdf}
\usepackage{amsmath}
\usepackage{color}
\usepackage{algorithm}
\usepackage{algpseudocode}

\usepackage{cite}
\newtheorem{problem}{Problem}

\newtheorem{remark}{Remark}
\newtheorem{theorem}{Theorem}

\def\BibTeX{{\rm B\kern-.05em{\sc i\kern-.025em b}\kern-.08em
    T\kern-.1667em\lower.7ex\hbox{E}\kern-.125emX}}
\usepackage{balance}
\begin{document}
\title{Finite-Horizon Discrete-Time Optimal Control for Nonlinear Systems under State and Control Constraints}
\author{Chuanzhi Lv, Hongdan Li, Huanshui Zhang
	
	\thanks{This work was supported by the 
		Original Exploratory Program Project of National Natural Science Foundation of China (62450004),  
		the Joint Funds of the National Natural Science Foundation of China (U23A20325),  
		the Major Basic Research of Natural Science Foundation of Shandong Province (ZR2021ZD14), 
		the National Natural Science Foundation of China (62103241),
		the Shandong Provincial Natural Science Foundation (ZR2021QF107), 
		and the Development Plan for Youth Innovation Teams in Higher Education Institutions in Shandong Province (2022KJ222). (\textit{Corresponding author: Huanshui Zhang.})
		
		C. Lv and H. Li are with the College of Electrical Engineering and Automation, Shandong University of Science and Technology, Qingdao 266590, China (e-mail: lczbyq@sdust.edu.cn; lhd2022@sdust.edu.cn).
		
		H. Zhang is with the College of Electrical Engineering and Automation, Shandong University of Science and Technology, Qingdao 266590, China, and  also with the School of Control Science and Engineering, Shandong University, Jinan, Shandong 250061, China (e-mail: hszhang@sdu.edu.cn).
		
		This work has been submitted to the IEEE for possible publication. Copyright may be transferred without notice, after which this version may no longer be accessible.
 	}
}

\maketitle

\begin{abstract}
	This paper addresses the optimal control problem of finite-horizon discrete-time nonlinear systems under state and control constraints.
	A novel numerical algorithm based on optimal control theory is proposed to achieve superior computational efficiency,  with the novelty lying in establishing a unified framework that integrates all aspects of algorithm design through the solution of forward and backward difference equations (FBDEs).
	Firstly, the state and control constraints are transformed using an augmented Lagrangian method (ALM), thereby decomposing the original optimal control problem into several optimization subproblems. 
	These subproblems are then reformulated as new optimal control problem, which are solved through the corresponding FBDEs, resulting in an algorithm with superlinear convergence rate. 
	Furthermore, the gradient and Hessian matrix are computed by iteratively solving FBDEs, thereby accelerating the optimization process.  
	The gradient is obtained through the standard Hamiltonian, while the Hessian matrix is derived by constructing a novel Hamiltonian specifically designed for second-order optimization, transforming each row into an iterative solution of a new set of FBDEs.
	Finally, the effectiveness of the algorithm is validated through simulation results in automatic guided vehicles (AGV) trajectory tracking control.
\end{abstract}


\section{Introduction}
\IEEEPARstart{O}{ptimal} control is fundamentally a type of dynamic optimization, which addresses problems involving time-dependent state variables and control inputs in a sequential or continuously evolving environment. Unlike static optimization, dynamic optimization accounts for both the immediate outcomes of decisions and their long-term impacts on future system states. 
Optimal control methods have found widespread application across diverse interdisciplinary domains, including engineering and artificial intelligence \cite{worthmann2015model, zhao2024survey, zhou2022swarm, wang2022geometrically, daoud2022simultaneous, li2018optimal}, with motion planning \cite{zhao2024survey}, trajectory optimization \cite{wang2022geometrically}, and path following \cite{daoud2022simultaneous} for AGV and aerial robots being notable examples that can be mathematically formulated as specific cases of such problems.
Due to the inherent limitations of physical systems and practical scenarios, dynamic systems are inevitably subject to state and control constraints, a subject  that has garnered significant attention from researchers. 

\par A classical numerical approach to solving optimal control problems is based on the Pontryagin's maximum principle, which provides first-order optimality condition and is also called the indirect method \cite{naidu2002optimal}.
The core of the method involves solving a set of FBDEs to obtain an exact solution.
However, the strong coupling between the forward and backward equations often presents significant challenges in the solution process. 
Furthermore, the maximum principle is limited in its ability to handle complex constraints.
Existing methods for handling constraints primarily modify the cost function. One approach involves constructing continuously integrable non-quadratic functionals to enforce the constraints \cite{lyshevski1998nonlinear,lyshevski2024analytic}, while another introduces penalty functions to penalize constraint violations \cite{ohtsuka1997real}.
At the algorithmic level, one can employ existing optimizers to solve a two-point boundary value problem and obtain the optimal solution.
However, this method is sensitive to the initial costate values and is typically applicable to affine nonlinear systems \cite{pagone2024continuous}. 
For more general systems, additional equilibrium condition within the optimality condition must be solved, which introduces a significant computational burden.
Additionally, a method of successive approximation (MSA) can be used to iteratively solve the FBDEs \cite{li2018maximum}. 
This method, essentially based on gradient descent, is feasible but tends to converge slowly.

\par A prominent approach is to reformulate the original problem as a nonlinear optimization problem, known as the direct method \cite{rao2009survey}. 
This approach offers advantages in terms of strong numerical stability and flexibility in handling complex constraints, and has been widely applied in practical engineering scenarios \cite{li2021design, cenerini2023model}. 
Commonly used optimization algorithms include sequential quadratic programming (SQP) and the interior-point method \cite{bazaraa2006nonlinear}.
Unlike general nonlinear optimization problems, optimal control problems are inherently dynamic.
Consequently, the gradient and Hessian matrix often need to be computed by existing optimizers during the solution process.
Two commonly used techniques for computing these derivatives are finite difference and automatic differentiation \cite{nocedal1999numerical}.
Notably, the finite difference method relies on the Taylor series expansion and is prone to truncation errors, which can reduce the accuracy of gradient and Hessian  matrix estimation. 
The automatic differentiation method, despite providing highly accurate derivatives, has its own limitations as it requires constructing a computational graph and depends on specialized optimization tools.

\par The combination of the strengths of the maximum principle and nonlinear optimization in a hybrid optimization method has attracted significant attention \cite{lin2014control, fu2021dynamic}.
By applying the maximum principle for the unified computation of the gradient of the cost function, this approach facilitates the use of standard gradient-based nonlinear optimization techniques to efficiently solve them \cite{teo2021applied}.
By replacing automatic differentiation or finite difference, this unified gradient computation not only ensures accuracy but also enhances computational efficiency, thereby accelerating the optimization process in gradient-based algorithms.
Relevant practical applications can be found in \cite{ren2023dynamic}.
As is well known, the second-order derivative information of the cost function can improve the convergence properties of the optimization process.
However, the current hybrid methods cannot be extended to second-order optimization algorithms, as the maximum principle does not provide a direct approach for computing the Hessian matrix, creating a significant obstacle.

\par The discussion above motivates us to develop efficient second-order algorithm for accelerating optimal control computation.
From the viewpoint of the maximum principle, the core of optimal control problems lies in solving FBDEs. In this paper, we frame the design of the optimization algorithm, as well as the computation of gradient and Hessian matrix, as the solution of FBDEs. 
In particular, the gradient and Hessian matrix computation does not require any external optimizers and can also be applied to other second-order algorithms.
This approach ensures a high level of uniformity in the algorithm design, making it stable, easy to implement, and capable of faster computation.
The main contributions of this paper are summarized as follows:
\begin{itemize}
	\setlength{\itemindent}{1.5mm}
	\item[1)]By applying an ALM, the original nonlinear optimal control problem is transformed into a series of optimization subproblems. These subproblems are then reformulated from the perspective of optimal control. The corresponding solutions are obtained by solving the FBDEs, and an optimization algorithm with a superlinear convergence rate is presented.
	\item[2)]The precise explicit formulas for the gradient and Hessian matrix, derived from calculus of variations, are proposed to accelerate the optimization process.
	The gradient is obtained through the standard Hamiltonian, while the Hessian matrix is derived by constructing a novel Hamiltonian specifically designed for second-order optimization, which transforms each row into an iterative solution of a new set of FBDEs.
	\item[3)]The proposed algorithm is applied to trajectory tracking control of AGV with state and control constraints. Simulation results demonstrate that this approach achieves superior computational efficiency.
\end{itemize}

\par The rest of this paper is organized as follows. In section II, the discrete-time nonlinear optimal control problem under state and control constraints is presented. 
In section III, the transcription of state and control constraints, the optimization algorithm design, and the formulas for computing the gradient and Hessian matrix are presented.
In Section IV, the proposed algorithm is applied to trajectory tracking control for AGV.
Finally, some concluding remarks are given in Section V.

\par \textit{Notations:} $\mathbb{R}^n$ denotes the $n$-dimensional Euclidean space; $A^{T}$ represents the transpose of matrix $A$; $\frac{\partial p}{\partial q}$ denotes partial derivative of function $p$ with respect to $q$. $\textbf{0}$ denotes the zero vector.


\section{Problem Statement}
Consider the following discrete-time nonlinear system:
\begin{equation}
	\begin{cases}
		x_{k+1} = f(x_k,u_k), \\
		x_0 = \xi,\ k=0,...,N,        \\
	\end{cases}                         \label{2_1}
\end{equation}
where $x_k \in \mathbb{R}^n$ is the state variable, $u_k \in \mathbb{R}^m$ is the control variable, $x_0$ is the given initial state, and $f\in\mathbb{R}^n$ is a given deterministic function. 
Define the following state and control constraints 
\begin{align}
	c_i(x_k,u_k)\leq 0, \ i=1,2,...,l.    \label{2_2}
\end{align}

\par The associated cost function is defined as follows:
\begin{align}
	J = \sum_{k=0}^{N}L(x_k,u_k),   \label{2_3}
\end{align}
where $L\in\mathbb{R}$ is the running cost. It is worth mentioning that the cost function (\ref{2_3}) considered here is in the Lagrange form. The Bolza form can be regarded as a special case of this formulation, obtained by setting $L(x_N,u_N)=\phi(x_N)$ at $k=N$. 
As observed from (\ref{2_1}) and (\ref{2_3}), the cost function $J$ can be expressed as a function of the initial state $x_0$ and the control sequence ($u_0,...,u_N$), i.e., $J=J(x_0,u_k)$.

\par Then the problem of this paper is formulated as follows:
\begin{problem}
	Find the optimal controller $u_k$ to minimize the cost function (\ref{2_3}) subject to the dynamical system (\ref{2_1}) and the state and control constraints (\ref{2_2}). \label{Problem1}
\end{problem}

\section{Numerical Algorithm based on Optimal Control Method}
In this section, the state and control constraints are firstly addressed through the use of ALM. A superlinearly convergent algorithm is then presented, accompanied by explicit iterative formulas for the computation of both the gradient and the Hessian matrix. Finally, a detailed description of the implementation process is provided.

\subsection{Constraints Transcription}
The ALM is an effective algorithm for solving constrained nonlinear optimization problems. 
Its fundamental idea is to incorporate the constraints directly into the objective function by introducing Lagrange multipliers along with additional penalty terms for any constraint violations. 
This approach transforms the original constrained problem into a series of unconstrained optimization problems \cite{bazaraa2006nonlinear}.
In order to address the Problem \ref{Problem1}, define the following augmented Lagrange function
\begin{align}
	\nonumber &\widetilde{J}(x_0,u_k,\gamma_{i,k},\rho_{i,k},\sigma) \\
	\nonumber = &J(x_0,u_k) + \sum_{i=1}^{l}\sum_{k=0}^{N}\gamma_{i,k}\big[c_i(x_k,u_k)+\rho_{i,k}^2\big]\\
	&\hspace*{3.96em}+ \frac{\sigma}{2}\sum_{i=1}^{l}\sum_{k=0}^{N}\big[c_i(x_k,u_k)+\rho_{i,k}^2\big]^2,   \label{3aa-1}
\end{align}
where $\gamma_{i,k}$ are the Lagrange multipliers, $\rho_{i,k}$ are the slack variables, and $\sigma > 0$ is the penalty parameter. 
To reduce the dimensionality of the optimization problem, $\rho_{i,k}$ is first minimized to obtain its optimal solution, which is then substituted into the function (4). 
Following this, the terms are reorganized to obtain the following reformulated optimization problem:
\begin{align}
	\hspace*{-0.02em}\nonumber \underset{u_k}{\text{min}} \  &\widetilde{J}(x_0,u_k,\gamma_{i,k},\sigma) =J(x_0,u_k) +\\
	&\frac{1}{2\sigma}\sum_{i=1}^{l}\sum_{k=0}^{N}\big\{\big[\text{max}(0,\gamma_{i,k}+\sigma c_i(x_k,u_k))\big]^2-\gamma_{i,k}^2\big\}.   \label{3aa-2}
\end{align}
The ALM transforms Problem \ref{Problem1} into a sequence of subproblems  (\ref{3aa-2}). 
By gradually increasing the penalty parameter $\sigma$ and sequentially solving these subproblems, the solutions of the unconstrained subproblems progressively approach the optimal solution of Problem \ref{Problem1}.
Note that in cases where high accuracy is not essential, a simpler exterior penalty method can replace the above approach.

\par The specific implementation steps of the ALM are presented below.

\par \textit{Step 1:} Let a penalty parameter $\sigma_1>0$, a termination scalar $\varepsilon>0$, a scalar $\beta>1$, and $j=1$. 
Initialize the control sequence $u_k^0,\ k=0,1,...,N$ and the multiplier sequence $\gamma_{i,k}^1,$ $i=1,...,l,\ k=0,1,...,N$.

\par \textit{Step 2:} Starting with the $u_k^{j-1}$, solve the following optimization problem: 
\begin{align}
	\underset{u_k}{\text{min}} \ \widetilde{J}(x_0,u_k,\gamma_{i,k}^j,\sigma_j),   \label{3a-2}  
\end{align}
where the optimal solution $u_k^{j}$ is obtained.

\par \textit{Step 3:} If $\sum_{i=1}^{l}\sum_{k=0}^{N}\big[\text{max}(c_i(x_k^j,u_k^j),-\frac{\gamma_{i,k}^j}{\sigma_j})\big]^2 <\varepsilon$, stop. Otherwise, go to Step 4.

\par \textit{Step 4:} Let $\gamma_{i,k}^{j+1}=\text{max}\big\{0,\gamma_{i,k}^j+\sigma_j c_i(x_k^j,u_k^j)\}$ and $\sigma_{j+1}=\beta\sigma_j$, replace $j$ by $j+1$, and go to Step 2.

\par The core of the four steps outlined above lies in solving optimization problem (\ref{3a-2}).
The design process of the corresponding algorithm will be presented below.

\subsection{Optimization Algorithm Design}
Let $\boldsymbol{u}$ denote the vector of control inputs over all time steps, i.e., $\boldsymbol{u}=[u_0^{T},...,u_N^{T}]^{T}$, and $\widetilde{J}(x_0,\boldsymbol{u},\gamma_{i,k}^j,\sigma_j) =\widetilde{J}(x_0,u_k,\gamma_{i,k}^j,\sigma_j)$.
The general form of the optimization algorithm is given by
\begin{align}
	\boldsymbol{u}_{t+1} = \boldsymbol{u}_t + s_t, \label{3b-1}
\end{align}
where $t=0,1,...,M-1$, and $s_t$ denotes the update step, with $t$ representing the current iteration of the optimization process and $M$ indicating the maximum number of iterations allowed. For example, in the case of the basic gradient descent method, the iteration step is represented by $-\varpi \nabla\widetilde{J}(x_0,\boldsymbol{u}_t,\gamma_{i,k}^j,\sigma_j)$, where $\varpi$ is a scalar positive real number. 

\par From the perspective of optimal control, the $\boldsymbol{u}_t$ and iteration step $s_t$ can be viewed as state and control variable, respectively, with the goal of minimizing the accumulated cost from step $t+1$ to step $M-1$, i.e., $\sum_{l=t+1}^{M-1}\widetilde{J}(x_0,\boldsymbol{u_l},\gamma_{i,k}^j,\sigma_j)$. At the same time, it is necessary to minimize the energy consumed by each iteration step. Accordingly, define the following cost function
\begin{align}
	\nonumber\widetilde{J}_M \hspace*{-0.1em}=\hspace*{-0.2em} \sum_{t=0}^{M-1}\big[\widetilde{J}(x_0,\boldsymbol{u}_t,\gamma_{i,k}^j,\sigma_j) &+\frac{1}{2}s_t^T\mathcal{R}s_t\big]
	\\
	&+\widetilde{J}(x_0,\boldsymbol{u}_M,\gamma_{i,k}^j,\sigma_j),               \label{3b-2}
\end{align}
where $\mathcal{R}$ is a positive definite matrix. 

\par The goal is to find the optimal control $s_t$ to minimize the cost function (\ref{3b-2}) subject to (\ref{3b-1}). By applying the maximum principle to the above problem, the costate equation and equilibrium condition can be derived as follows:
\begin{align}
	\zeta_t &= \zeta_{t+1} + \nabla \widetilde{J}(x_0,\boldsymbol{u_t},\gamma_{i,k}^j,\sigma_j),      \label{3b-3}
	\\
	0&=\mathcal{R}s_t + \zeta_t,                \label{3b-4}
\end{align}
where $\zeta_t$ is costate variable.
The equations (\ref{3b-1}) and (\ref{3b-3})-(\ref{3b-4}) together form the FBDEs, which exhibit strong coupling. 
It is evident from the FBDEs that when the state $\boldsymbol{u}_t$ in the system (\ref{3b-1}) stabilizes, the terminal state $\boldsymbol{u}_M$ represents the optimal solution to the problem (\ref{3a-2}). 
In general, except for linear systems, it is almost impossible to obtain an analytical solution. In practical computations, the MSA \cite{xu2023optimization} can be employed for numerical solutions. However, it should be noted that the MSA is essentially a gradient descent method, which results in a relatively slow convergence rate. 
To enhance computational efficiency, the following iterative algorithm with superlinear convergence rate can be obtained by combining the Taylor expansion:
\begin{align}
	\boldsymbol{u}_{t+1} &= \boldsymbol{u}_t - h_t, \label{3b-5}
	\\
	h_m &=(\mathcal{R} + \nabla^2 \widetilde{J}(t))^{-1}(\nabla \widetilde{J}(t) + \mathcal{R}h_{m-1}), \label{3b-6}
	\\
	h_0 &= (\mathcal{R} + \nabla^2 \widetilde{J}(t))^{-1}\nabla \widetilde{J}(t), \label{3b-7}
\end{align}
where $\nabla^2 \widetilde{J}(t)=\nabla^2 \widetilde{J}(x_0,\boldsymbol{u_t},\gamma_{i,k}^j,\sigma_j)$, $\nabla \widetilde{J}(t)=\nabla \widetilde{J}(x_0, \boldsymbol{u_t},$ $\gamma_{i,k}^j,\sigma_j)$, $m=1,...,t$, and $t=0,1,...,M-1$. 
For the detailed derivation process and analysis, refer to \cite{zhang2024optimization}.

\par As shown in the formulas (\ref{3b-5})-(\ref{3b-7}), each iteration step $h_t$ requires $t+1$ inner loop iterations to obtain the result.
The key to the algorithm is to compute the gradient $\nabla \widetilde{J}(t)$ and Hessian matrix $\nabla^2 \widetilde{J}(t)$ of the cost function $\widetilde{J}(x_0,\boldsymbol{u}_t,\gamma_{i,k}^j,\sigma_j)$ at each iteration step.
The explicit iterative results of computing the gradient and Hessian matrix for dynamic optimization will be presented below.

\subsection{Computation of Gradient and Hessian Matrix}
Unlike static optimization, the cost function $\widetilde{J}(x_0,u_k,\gamma_{i,k}^j,$ $\sigma_j)$ in the optimization problem (\ref{3a-2}) involves a dynamical system. Common methods, such as automatic differentiation or finite difference, are typically used to compute its gradient and Hessian matrix. 
However, these techniques make the computation with respect to the control variable $u_k$ generally more complex and time-consuming.
The maximum principle addresses this challenge by introducing Lagrange multipliers \cite{naidu2002optimal}, which transform the gradient problem into solving a set of FBDEs, thereby significantly improving computational efficiency. In the following, the explicit expression for the gradient is derived for the case of a one-dimensional controller, with natural extensions to higher-dimensional cases.
For clarity and readability, the cost function $\widetilde{J}(x_0,u_k,\gamma_{i,k}^j,\sigma_j)$ in (\ref{3a-2}) is denoted as $\widetilde{J}(x_0,u_k,\gamma_{i,k},\sigma)$ in the following derivation of the gradient and Hessian matrix.

\begin{theorem} \label{theorem1}
	For the optimization problem (\ref{3a-2}), the gradient of $\widetilde{J}(x_0,u_k,\gamma_{i,k},\sigma)$ with respect to $u_k$ is given by 
	\begin{align}
		\nabla\widetilde{J}(x_0,u_k,\gamma_{i,k},\sigma) = \left[\frac{\partial \mathcal{H}(0)}{\partial u_0},...,\frac{\partial \mathcal{H}(N)}{\partial u_{N}}\right]^T,    \label{3c-1}
	\end{align}
	where Hamiltonian $\mathcal{H}(k)$, for $k=0,1,...,N$, is given by 
	\begin{align}
		\nonumber \mathcal{H}(k) &= \mathcal{H}(x_k,u_k,\lambda_{k+1},\gamma_{i,k},\sigma)
		\\
		\nonumber&= L(x_k,u_k) + \lambda_{k+1}^{T}f(x_k,u_k) + 
		\\ &\hspace*{1.0em}\frac{1}{2\sigma}\sum_{i=1}^{l}\big\{\big[\text{max}(0,\gamma_{i,k}+\sigma c_i(x_k,u_k))\big]^2 - (\gamma_{i,k})^2\big\}. \label{3c-4}
 	\end{align}
	Moreover, the costate $\lambda_k$, for $k=1,...,N+1$,  satisfies
	\begin{align}
		\lambda_k &= \frac{\partial \mathcal{H}(x_k,u_k,\lambda_{k+1},\gamma_{i,k},\sigma)}{\partial x_k}, \label{3c-5}
	\end{align}
	with the terminal value $\lambda_{N+1} = \boldsymbol{0}$.
\end{theorem}

\par \textit{Proof} 1:	The computation of gradient (\ref{3c-1}) can be derived by applying the maximum principle. Specifically, by introducing Lagrange multiplier $\lambda_k$ in (\ref{3c-4}) to the dynamical system (\ref{2_1}), the cost function $\widetilde{J}(x_0,u_k,\gamma_{i,k},\sigma)$ can be reformulated in the following augmented form:
	\begin{align}
		\nonumber \hspace*{-0.7em}&\widetilde{J}(x_0,u_k,\lambda_{k+1},\gamma_{i,k},\sigma) 
		\\
		\nonumber \hspace*{-0.7em}=& \sum_{k=0}^{N}\left[L(x_k,u_k) + \lambda_{k+1}^T(f(x_k,u_k)-x_{k+1})\right]  +
		\\
		\hspace*{-0.7em}& \frac{1}{2\sigma}\sum_{i=1}^{l} \sum_{k=0}^{N} \big\{\big[\text{max}(0,\gamma_{i,k}+\sigma c_i(x_k,u_k))\big]^2 - (\gamma_{i,k})^2\big\}. \label{3c-6}
	\end{align}
	The derivation follows the principle of the calculus of variations, where the necessary conditions for optimality are obtained by analyzing the variations of the augmented cost function (\ref{3c-6}), ultimately leading to the explicit expressions of the gradient (\ref{3c-1})-(\ref{3c-5}).
	The proof is complete.

\par The essence of computing the gradient $\nabla \widetilde{J}(x_0,u_k,\gamma_{i,k},\sigma)$ lies in solving the optimal control problem associated with the cost function $\widetilde{J}$ subject to the forward difference equation (\ref{2_1}). 
Note that the maximum principle does not directly provide a method for computing the Hessian matrix. 
As seen from the gradient expressions (\ref{3c-1})-(\ref{3c-5}), when computing the Hessian matrix, the elements of the gradient (\ref{3c-1}) are constrained not only by the forward difference equation (\ref{2_1}) but also by the backward difference equation (\ref{3c-5}). Therefore, the result of Theorem \ref{theorem1} cannot be directly extended.

\par In the following, explicit iterative formulas for computing the Hessian matrix of the cost function $\widetilde{J}(x_0,u_k,\gamma_{i,k},\sigma)$ are derived using the calculus of variations.

\begin{theorem}
	For the optimization problem (\ref{3a-2}), the Hessian matrix of $\widetilde{J}(x_0,u_k,\gamma_{i,k},\sigma)$ with respect to $u_k$ is given by 
	\begin{align}
		\nabla^2\widetilde{J}(x_0,u_k,\gamma_{i,k},\sigma) &= \renewcommand{\arraystretch}{1.43} 
		\setlength{\arraycolsep}{3pt}
		\left[
		\begin{array}{cccc}
			\frac{\partial \mathcal{H}_0(0)}{\partial u_0} & \frac{\partial \mathcal{H}_0(1)}{\partial u_1} & \cdots & \frac{\partial \mathcal{H}_0(N)}{\partial u_N} \\
			\frac{\partial \mathcal{H}_1(0)}{\partial u_0} & \frac{\partial \mathcal{H}_1(1)}{\partial u_1} & \cdots & \frac{\partial \mathcal{H}_1(N)}{\partial u_N} \\
			\vdots & \vdots & \ddots & \vdots \\
			\frac{\partial \mathcal{H}_N(0)}{\partial u_0} & \frac{\partial \mathcal{H}_N(1)}{\partial u_1} & \cdots & \frac{\partial \mathcal{H}_N(N)}{\partial u_N}
		\end{array}
		\right]
		\renewcommand{\arraystretch}{1} 
		\setlength{\arraycolsep}{5pt},  \label{3c-7}
	\end{align}
	where Hamiltonian $\mathcal{H}_j(k)$, for $j=0,1,...,N$, $k=0,1,...,N$, is given by 
	\begin{align}
		\nonumber \hspace*{-0.66em}\mathcal{H}_j(k) &= \mathcal{H}_j(x_k,u_k,\lambda_{k+1},\gamma_{i,k},\sigma,\eta_{k+1},\mu_k)  
		\\
		\nonumber &= \Phi(x_j,u_j,\lambda_{j+1},\gamma_{i,j},\sigma)\boldsymbol{I}_{\left\{j=k\right\}}  
		\\
		&\hspace*{1em}+ \eta_{k+1}^{T}f(x_k,u_k)+ \mu_k^T\Psi(x_k,u_k,\lambda_{k+1},\gamma_{i,k},\sigma),     \label{3c-10}
	\end{align}
	with 
	\begin{align}
		\Phi(x_j,u_j,\lambda_{j+1},\gamma_{i,j},\sigma) &= \frac{\partial \mathcal{H}(x_j,u_j,\lambda_{j+1},\gamma_{i,j},\sigma)}{\partial u_j},  \label{3c-11}
		\\
		\Psi(x_k,u_k,\lambda_{k+1},\gamma_{i,k},\sigma) &= \frac{\partial \mathcal{H}(x_k,u_k,\lambda_{k+1},\gamma_{i,k},\sigma)}{\partial x_k}.   \label{3c-12}
	\end{align}
	Moreover, the costates $\mu_k$, for $k=0,1,...,N$, and $\eta_k$, for $k=1,...,N+1$ satisfy 
	\begin{align}
		\mu_{k+1} &= \frac{\partial \mathcal{H}_j(x_k,u_k,\lambda_{k+1},\gamma_{i,k},\sigma,\eta_{k+1},\mu_k)}{\partial \lambda_{k+1}}, \label{3c-13}
		\\
		\eta_k &= \frac{\partial \mathcal{H}_j(x_k,u_k,\lambda_{k+1},\gamma_{i,k},\sigma,\eta_{k+1},\mu_k)}{\partial x_k}, \label{3c-14}
	\end{align}
	with the initial value $\mu_0=\boldsymbol{0}$ and the terminal value $\eta_{N+1}=\boldsymbol{0}$.
\end{theorem}

\par \textit{Proof} 2:	To compute the Hessian matrix $\nabla^2\widetilde{J}(x_0,u_k,\gamma_{i,k},$ $\sigma)$ in (\ref{3c-7}), we proceed row by row, where each row corresponds to the gradient of the cost function (\ref{3c-11}) with respect to the control variable $u_k$.
	
	\par Consider the $j$-th row of the Hessian matrix. Define the cost function of the $j$-th row as
	\begin{align}
		V(\boldsymbol{u}) = \sum_{k=0}^{N}\Phi(x_k,u_k,\lambda_{k+1},\gamma_{i,k},\sigma), \label{3c-15}
	\end{align}
	where $\Phi(x_k,u_k,\lambda_{k+1},\gamma_{i,k},\sigma)=0$, $j\neq k$. By perturbing $\boldsymbol{u}$ in the cost function (\ref{3c-15}) with $\epsilon\cdot \delta\boldsymbol{u}$, where $\epsilon>0$ is a small real number, and $\delta\boldsymbol{u}=[\delta u_0, \delta u_1,...,\delta u_N]^{T}$, we can obtain
	\begin{align}
		\nonumber \boldsymbol{u_\epsilon} = \boldsymbol{u} + \epsilon\cdot \delta\boldsymbol{u},
	\end{align}
	where $\boldsymbol{u_\epsilon}=[u_{0,\epsilon},u_{1,\epsilon},...,u_{N,\epsilon}]^T$. Similarly, the dynamical system (\ref{2_1}) , the costate equation (\ref{3c-5}), and the cost function (\ref{3c-15}) are also perturbed, respectively, resulting in the following:
	\begin{align}
		\nonumber x_{k+1,\epsilon} &= f(x_{k,\epsilon},u_{k,\epsilon}),
		\\
		\nonumber \lambda_{k,\epsilon} &= \Psi(x_{k,\epsilon},u_{k,\epsilon},\lambda_{k+1,\epsilon},\gamma_{i,k},\sigma),
		\\
		V(\boldsymbol{u_\epsilon}) &= \sum_{k=0}^{N}\Phi(x_{k,\epsilon},u_{k,\epsilon},\lambda_{k+1,\epsilon},\gamma_{i,k},\sigma).  \label{3c-16}
	\end{align}
	Furthermore, the variation of the state $\delta x_k$ is given by 
	\begin{align}
		\nonumber \delta x_{k+1} &= \frac{\partial x_{k+1,\epsilon}}{\partial\epsilon}\Big|_{\epsilon=0}
		\\
		&= \frac{\partial f(x_k,u_k)}{\partial x_k}\delta x_k + \frac{\partial f(x_k,u_k)}{\partial u_k}\delta u_k,   \label{3c-17}
	\end{align}
	and the variation of the costate $\delta \lambda_k$ is given by
	\begin{align}
		\nonumber \hspace*{-0.5em}\delta \lambda_k = \frac{\partial \lambda_{k,\epsilon}}{\partial\epsilon}\Big|_{\epsilon=0} &= \frac{\partial \Psi(x_k,u_k,\lambda_{k+1},\gamma_{i,k},\sigma)}{\partial x_k}\delta x_k  
		\\
		\nonumber &\hspace*{1.3em}+\frac{\partial \Psi(x_k,u_k,\lambda_{k+1},\gamma_{i,k},\sigma)}{\partial u_k}\delta u_k 
		\\
		&\hspace*{1.3em} +\frac{\partial \Psi(x_k,u_k,\lambda_{k+1},\gamma_{i,k},\sigma)}{\partial \lambda_{k+1}}\delta \lambda_{k+1}.          \label{3c-18}
	\end{align}
	Combining the cost function (\ref{3c-15}) and (\ref{3c-16}), we can obtain
	\begin{align}
		\nonumber \hspace*{-1.0em}(\nabla V(\boldsymbol{u}))^T\delta\boldsymbol{u}&= \lim\limits_{\epsilon \to 0} \frac{V(\boldsymbol{u_\epsilon})-V(\boldsymbol{u})}{\epsilon}
		\\
		\nonumber &=\sum_{k=0}^{N}\big[\frac{\partial \Phi(x_k,u_k,\lambda_{k+1},\gamma_{i,k},\sigma)}{\partial x_k}\delta x_k 
		\\
		\nonumber &\hspace*{2.7em} +\frac{\partial \Phi(x_k,u_k,\lambda_{k+1},\gamma_{i,k},\sigma)}{\partial u_k}\delta u_k 
		\\
		&\hspace*{2.7em} +\frac{\partial \Phi(x_k,u_k,\lambda_{k+1},\gamma_{i,k},\sigma)}{\partial \lambda_{k+1}}\delta \lambda_{k+1}\big].   \label{3c-19}
	\end{align}
	Defining the Hamiltonian as is shown in (\ref{3c-10}) and substituting (\ref{3c-10}) into (\ref{3c-19}) results in:
	\begin{align}
		\nonumber \hspace*{-0.0em}&(\nabla V(\boldsymbol{u}))^T\delta\boldsymbol{u} 
		\\
		\nonumber =& \sum_{k=0}^{N}\big[\frac{\partial \mathcal{H}_j(k)}{\partial x_{k}}\delta x_k - \eta_{k+1}^T\frac{f(x_k,u_k)}{\partial x_{k}}\delta x_k 
		\\
		\nonumber &\hspace*{0.0em}- \mu_{k}^T\frac{\Psi(x_k,u_k,\lambda_{k+1},\gamma_{i,k},\sigma)}{\partial x_{k}}\delta x_k + \frac{\partial \mathcal{H}_j(k)}{\partial u_{k}}\delta u_k 
		\\
		\nonumber &\hspace*{0.0em}-\eta_{k+1}^T\frac{f(x_k,u_k)}{\partial u_{k}}\delta u_k - \mu_{k}^T\frac{\Psi(x_k,u_k,\lambda_{k+1},\gamma_{i,k},\sigma)}{\partial u_{k}}\delta u_k 
		\\
		&\hspace*{-0.0em}+\frac{\partial \mathcal{H}_j(k)}{\partial \lambda_{k+1}}\delta \lambda_{k+1} - \mu_{k}^T\frac{\Psi(x_k,u_k,\lambda_{k+1},\gamma_{i,k},\sigma)}{\partial \lambda_{k}}\delta \lambda_{k+1}\big]. \label{3c-20}
	\end{align}
	Defining the costates $\mu_k$ and $\eta_k$ as given in (\ref{3c-13}) and (\ref{3c-14}), and combining the variations (\ref{3c-17}) and (\ref{3c-18}), (\ref{3c-20}) can be rearranged as:
	\begin{align}
		\nonumber (\nabla V(\boldsymbol{u}))^T\delta\boldsymbol{u} &= \sum_{k=0}^{N}\frac{\partial \mathcal{H}_j(k)}{\partial u_{k}}\delta u_k + (\eta_0\delta x_0-\eta_{N+1}\delta x_{N+1})
		\\
		\nonumber &\hspace*{3.0em}+ (\mu_{N+1}\delta \lambda_{N+1}-\mu_0\delta \lambda_0)
		\\
		\nonumber &=\sum_{k=0}^{N}\frac{\partial \mathcal{H}_j(k)}{\partial u_{k}}\delta u_k.
		\\
		\nonumber &=\left[\frac{\partial \mathcal{H}_j(0)}{\partial u_{0}},...,\frac{\partial \mathcal{H}_j(N)}{\partial u_{N}}\right]\delta\boldsymbol{u}.
	\end{align}
	Considering the arbitrariness of $\delta \boldsymbol{u}$, the gradient of the cost function (\ref{3c-11}) with respect to the control variable $u_k$ is 
	\begin{align}
		\nonumber \nabla V(\boldsymbol{u}) = \left[\frac{\partial \mathcal{H}_j(0)}{\partial u_{0}},...,\frac{\partial \mathcal{H}_j(N)}{\partial u_{N}}\right]^T.
	\end{align}
	Therefore, the expression for the Hessian matrix in (\ref{3c-7}) holds.
	The proof is complete.

\begin{remark}
	The computational complexity of the gradient and Hessian matrix primarily comes from iteratively solving the FBDEs. 
	For gradient computation, assuming that only matrix addition and multiplication are involved, and that the state dimension $n$ is greater than the control dimension $m$, the computational complexities of both the forward and backward iterations in (\ref{2_1}) and (\ref{3c-5}) are the same. 
	Specifically, for each time step, the complexity of both iterations is $O(n^2)$. Therefore, the total complexity for gradient computation is $O(Nn^2)$, where $N$ is the total sampling time.
	For computing the Hessian matrix, each row involves iteratively sovling equations (\ref{3c-13}) and (\ref{3c-14}). 
	Since these equations are similar to those in gradient computation, the complexity analysis follows the same pattern. 
	As a result, the computational complexity for the Hessian matrix is $O(N^2n^2)$.
\end{remark}

\begin{algorithm}[tb]
	\renewcommand{\thealgorithm}{}
	\caption{\hspace{-0.6em} \textbf{1} Numerical algorithm with constraints.} 
	\label{alg:Framwork}
	\begin{algorithmic}[1]
		\Require 
		initial penalty factor $\sigma_1$, penalty factor update constant $\beta$, stopping criterion $\varepsilon$, $\bar \varepsilon$, initial state $x_{0}$, total sampling time $N$, prediction horizon $N_p$, initial value of   $\boldsymbol{u}_0$ at each sampling time;
		\Ensure  
		Optimal controller $u_k$ for each sampling time $k=0,...,N$ ;
		\State Perform the following optimization algorithm at each sampling time $k$.
		\For{$j=1$ \textbf{to} $\#step$}
		\For{$t=0$ \textbf{to} $M-1$}
		\State Compute the forward difference equation (\ref{2_1});
		\State Compute the backward difference equation (\ref{3c-5});
		\State Compute the gradient $\nabla \widetilde{J}(t)\hspace*{-0.26em}=\hspace*{-0.26em}\nabla\widetilde{J}(x_0,\boldsymbol{u}_t,\gamma_{i,k}^j,\sigma_j)$   
		\hspace*{3em}in (\ref{3c-1})-(\ref{3c-5}); 
		\If {$||\nabla\widetilde{J}(t)||^2<\bar \varepsilon$ }
		\State Break;
		\EndIf
		\State Compute the forward difference equation (\ref{3c-13});
		\State Compute the backward difference equation (\ref{3c-14});
		\State Compute the Hessian matrix $\nabla^2\widetilde{J}(t)=\nabla^2\widetilde{J}(x_0,$ \hspace*{3em}$\boldsymbol{u}_t,\gamma_{i,k}^j,\sigma_j)$ in (\ref{3c-7})-(\ref{3c-14});
		\State $h_0 = (\mathcal{R}+\nabla^2\widetilde{J}(t))^{-1}\nabla \widetilde{J}(t)$;
		\For{$m=1$ \textbf{to} $i$}
		\State $h_m = (\mathcal{R}+\nabla^2\widetilde{J}(t))^{-1}[\nabla \widetilde{J}(t) + \mathcal{R}h_{m-1}]$;
		\EndFor
		\State $\boldsymbol{u}_{t+1} = \boldsymbol{u}_t - h_t$;
		\EndFor
		\If {$\sum_{i=1}^{l}\sum_{k=0}^{N}\big[\text{max}(c_i(x_k^j,u_k^j),-\frac{\gamma_{i,k}^j}{\sigma_j})\big]^2 <\varepsilon$ }
		\State Break;
		\Else 
		\State $\gamma_{i,k}^{j+1} = \text{max}(0,\gamma_{i,k}^j + \sigma_jc_i(x_k^j,u_k^j))$;
		\State $\sigma_{j+1} = \beta\sigma_j$;
		\EndIf
		\EndFor
		\State Extract the first controller $u_k$ from $\boldsymbol{u}_{t+1}$;
		\State Update the state $x_k$ for the next sampling time;
	\end{algorithmic}
\end{algorithm}

\subsection{Computational Implementation}
The above describes an algorithm for solving optimal control over a finite time horizon. However, in practical implementation, due to error accumulation and model misspecification, solving the optimal control problem over a large time horizon may not be feasible \cite{weinan2022empowering}. Therefore, the concept of model predictive control (MPC) can be applied to address this issue. 
In MPC, at each sampling time $k=0,...,N$, the following optimal control problem, which builds upon (\ref{2_1})-(\ref{2_3}), is solved over a short prediction horizon $s=0,1,...,N_p$, where $N_p$ denotes the length of the prediction horizon:
\begin{alignat}{1}
	&\underset{u_{s|k}}{\text{min}} \  J_{N_p} = \sum_{s=0}^{N_p}L(x_{s|k},u_{s|k}),                   \label{3d-1}        
	\\
	&\hspace*{0.01em}  \mbox{s.t.}\ \hspace*{0.30em} x_{s+1|k} = f(x_{s|k}, u_{s|k}),                                            \label{3d-2}
	\\
	\nonumber & \hspace*{1.83em} x_{0|k} = \xi,\ k = 0,1,...,N,\ s = 0,1,...,N_p,   
	\\
	& \hspace*{1.83em}  c_i(x_k,u_k)\leq 0, \ i=1,2,...,l.     \label{3d-3}
\end{alignat}
The resulting optimal control sequence $\{ u^*_{s|k} \}_{s=0}^{N_p}$ is obtained. In practice, only the first control input $u^*_{0|k}$ is applied to the system, and the time horizon is shifted forward to the next time step $k+1$, where a new optimization problem is solved.
The details are summarized in Algorithm 1.

\begin{remark}
	The explicit iterative formulas for both the gradient and Hessian matrix of the cost function (\ref{2_3}) are obtained through the framework of FBDEs. By defining the standard Hamiltonian (\ref{3c-4}), the gradient computation is transformed into an iterative solution of FBDEs (\ref{2_1}) and (\ref{3c-5}). Based on this gradient computation, a novel Hamiltonian (\ref{3c-10}) is constructed to perform a second-order optimization of (\ref{3c-1}), converting each row of the Hessian matrix (\ref{3c-7}) into the iterative solution of a new set of FBDEs (\ref{3c-13}) and (\ref{3c-14}).
	It is worth noting that the forward difference equation (\ref{3c-13}) is derived from the backward difference equation (\ref{3c-5}), while the backward difference equation (\ref{3c-14}) is derived from the forward difference equation (\ref{2_1}).
\end{remark}
\begin{remark}
	The algorithm's high efficiency can be attributed to two key factors. Firstly, the algorithm (\ref{3b-5})-(\ref{3b-7}) exhibits superlinear convergence rate \cite{zhang2024optimization}, allowing for rapid convergence to the optimal solution, while maintaining a simple and efficient code implementation. 
	Secondly, the explicit iterative formulas for the gradient and Hessian matrix of the cost function $\widetilde{J}(x_0,u_k,\gamma_{i,k}^j,\sigma_j)$ in (\ref{3a-2}) can accelerate the computation process, eliminate the need for optimization tools such as automatic differentiation or finite difference, and reduce both computational and storage requirements, thereby improving overall efficiency.
\end{remark}

\section{Trajectory Tracking Control for AGV}
\subsection{Description of AGV system}
\begin{figure}[t]
	\centering
	\includegraphics[width=0.23\textwidth]{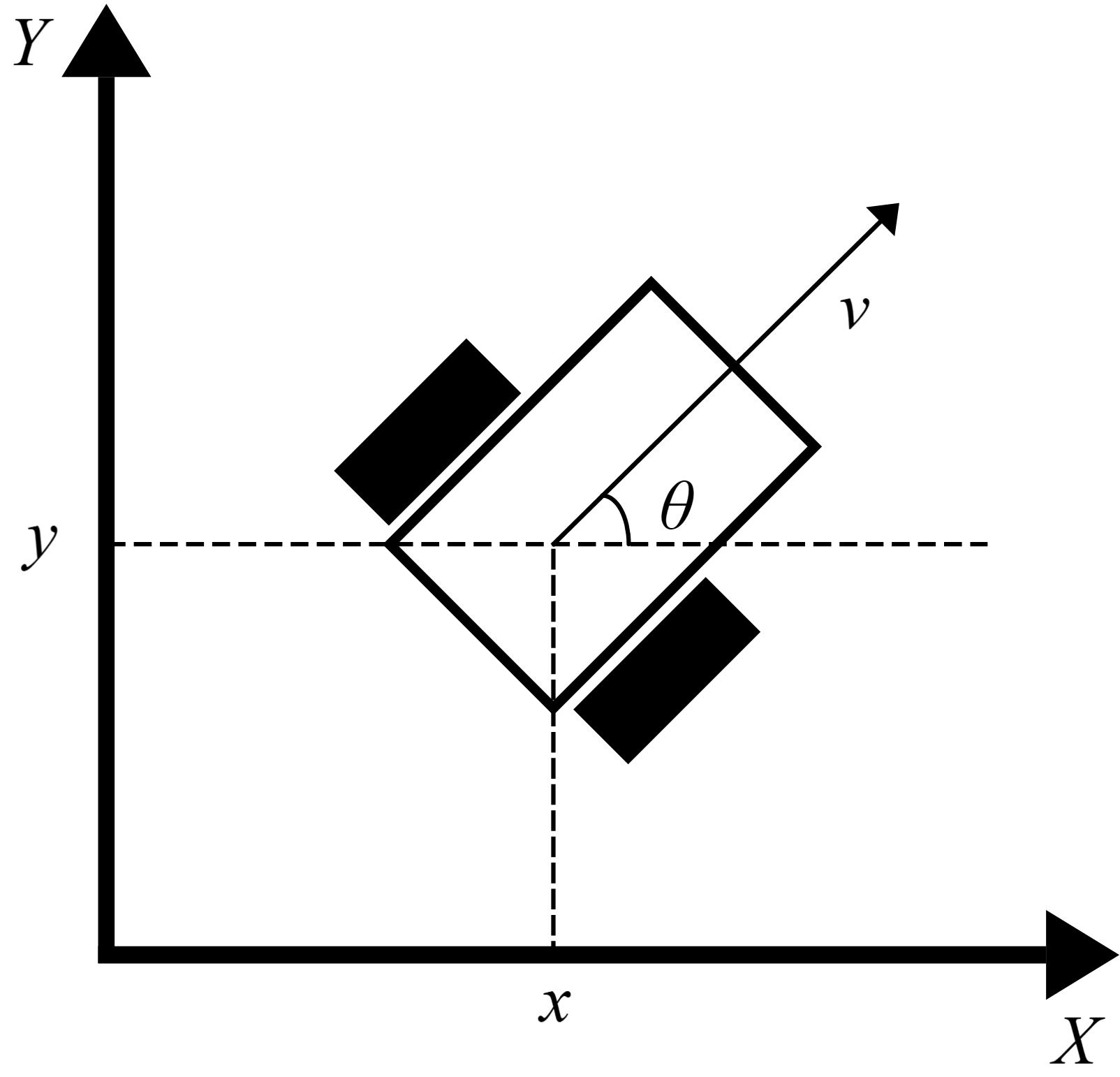}
	\captionsetup{font={footnotesize}}
	\caption{Model schematic.} \label{fig4a-2}
\end{figure}

The kinematic model of the AGV is given by
\begin{equation}
	\begin{cases}
		\dot x_t = v_t{\rm cos}\theta_t, \\
		\dot y_t = v_t{\rm sin}\theta_t,\\
		\dot \theta_t = \omega_t,
	\end{cases}                         \label{4a-1}
\end{equation}
where $t\in[0,t_f]$, $[x_t,y_t]^T$ is the position of the AGV, $\theta_t$ is the orientation, $[x_t,y_t,\theta_t]^T$ denotes the pose of the AGV with respect to the Cartesian coordinate system, with the initial pose $[x_0,y_0,\theta_0]^T$ known; $v_t$ and $\omega_t$ are the linear and angular velocities, respectively. 
The model schematic of the AGV is shown in Fig. \ref{fig4a-2}.

\begin{table}[t]
	\centering
	\captionsetup{font={footnotesize}}
	\caption{Comparison of total computation time and average computation time for Algorithm 1, fmincon, and Casadi in the MPC implementation for the trajectory tracking problem (\ref{4a-2})-(\ref{4a-3}).}
	\begin{tabular}{|c|c|c|c|c|c|c|c|c|c|c|c|}
		\hline
		\multicolumn{2}{|c|}{\multirow{1}{*}{Algorithm}} &\multirow{1}{*}{Total  time [s]}&\multirow{1}{*}{Average time [s]} \\
		\hline
		\multirow{2}{*}{fmincon} & Interior Point & 3.1836 & 0.0199 \\
		\cline{2-4}	
		\multirow{2}{*}{} & SQP & 1.2997 & 0.0081 \\
		\hline
		
		\multirow{1}{*}{Casadi} & IPOPT & 1.3768 & 0.0086 \\
		\hline
		
		\multicolumn{2}{|c|}{\multirow{1}{*}{Algorithm 1}} &\multirow{1}{*}{0.2809}&\multirow{1}{*}{0.0018} \\
		
		\hline
	\end{tabular}
	\label{table4a-1}
\end{table}

\begin{figure}[t]
	\centering
	\subfigure[\hspace*{0.0em}]{
		\includegraphics[width=0.45\textwidth]{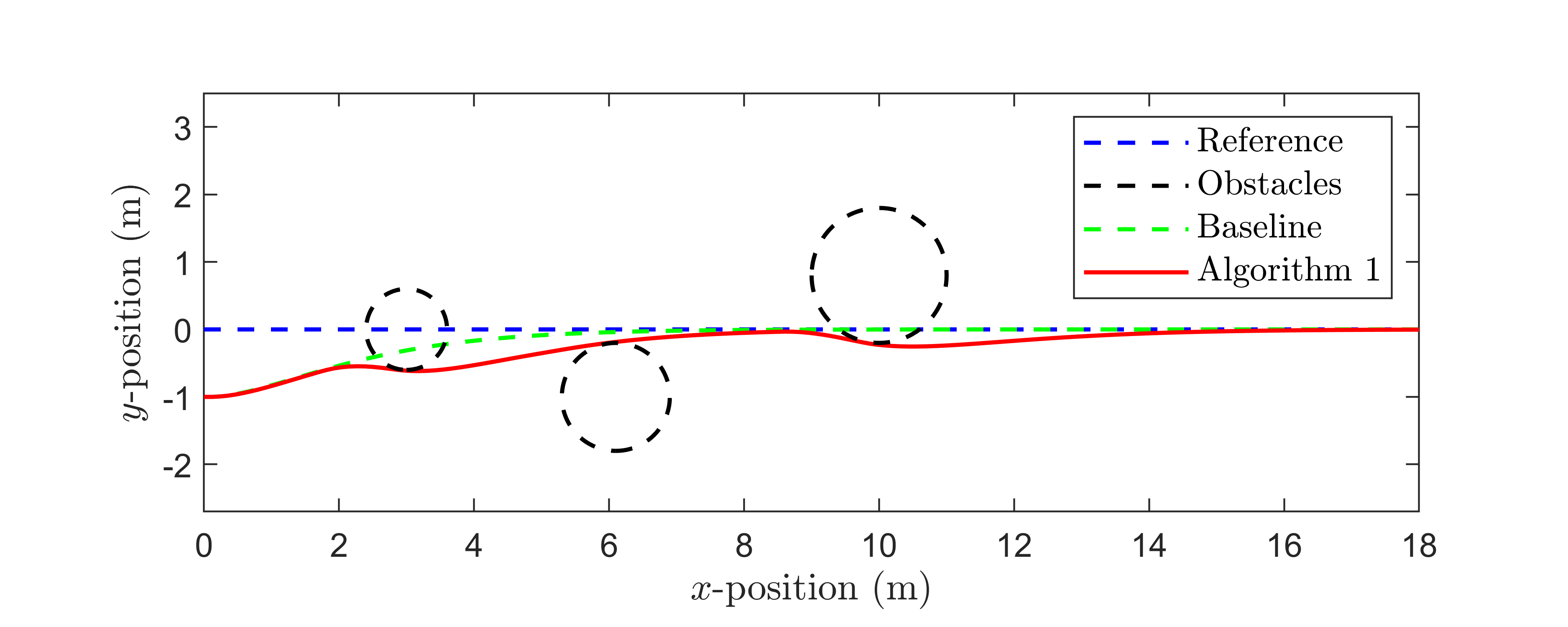}
	}
	\subfigure[\hspace*{0.0em}]{
		\includegraphics[width=0.45\textwidth]{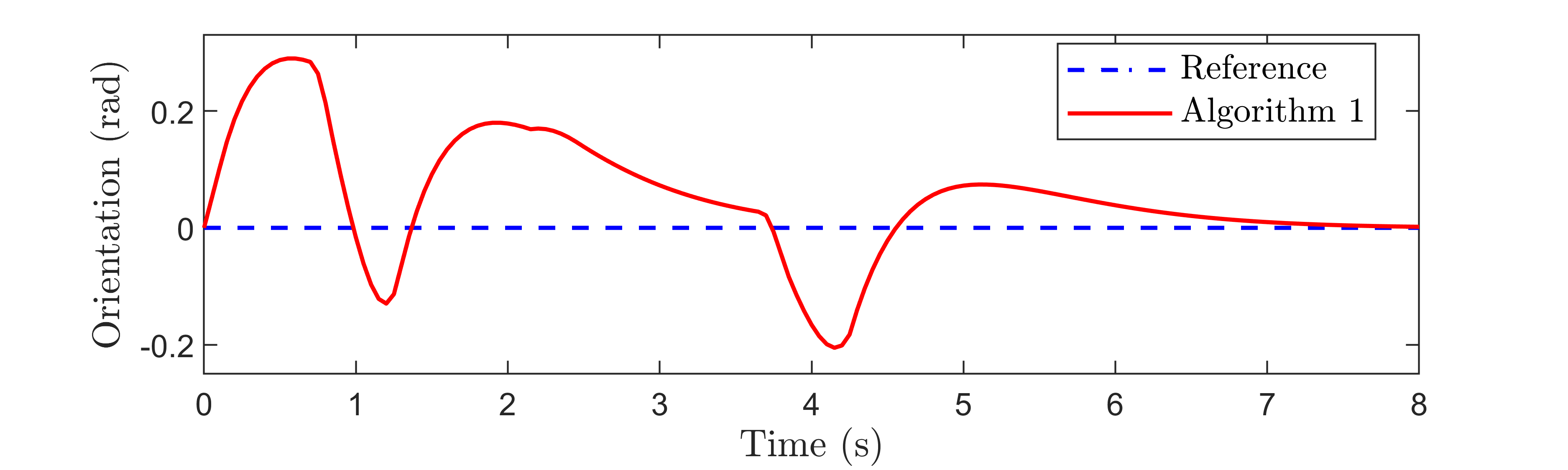}
	}
	\subfigure[\hspace*{0.0em}]{
		\includegraphics[width=0.45\textwidth]{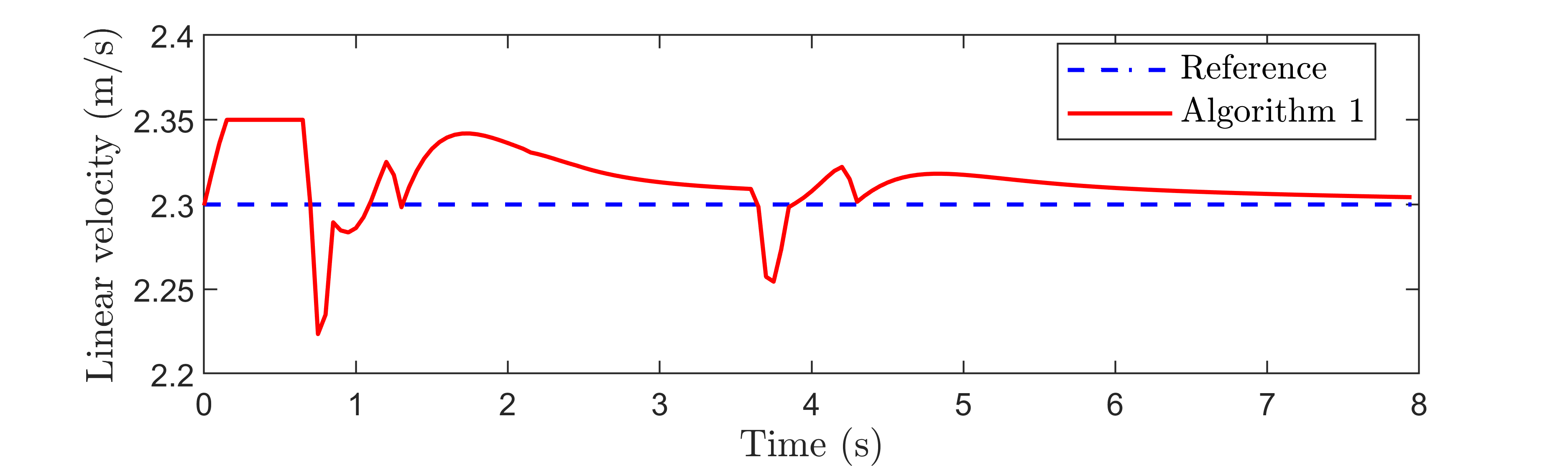}
	}
	\subfigure[\hspace*{0.0em}]{
		\includegraphics[width=0.45\textwidth]{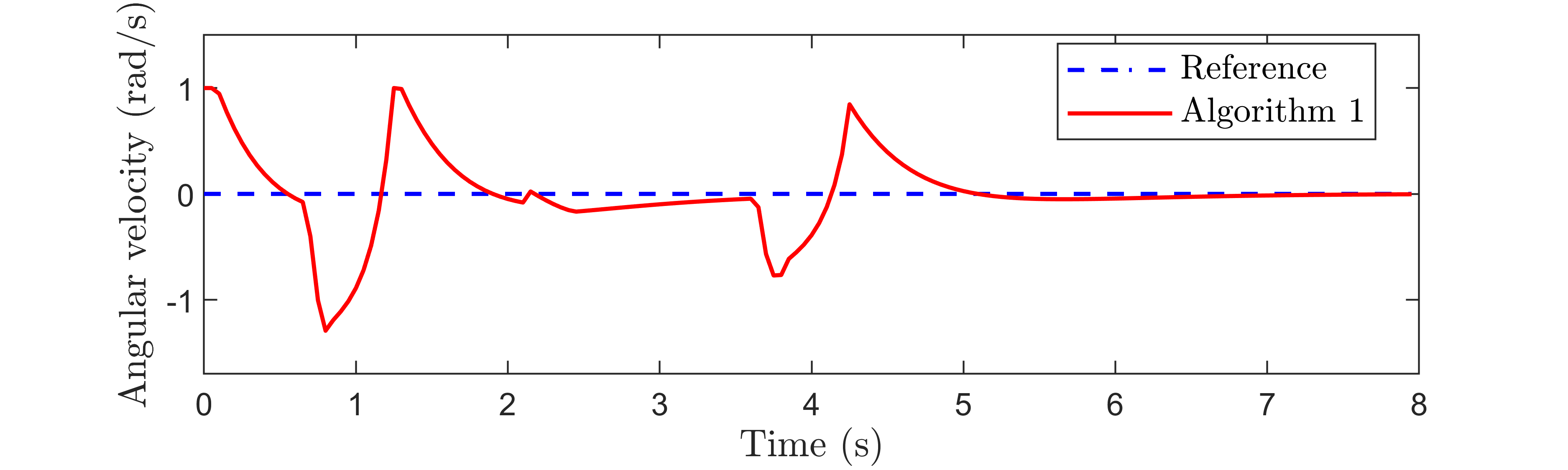}
	}
	\captionsetup{font={footnotesize}}
	\caption{Control performance of Algorithm 1 for the trajectory tracking problem (\ref{4a-2})-(\ref{4a-3}). (a) Tracking results. (b) Orientation $\omega_k$. (c) Linear velocity $v_k$. (d) Angular velocity $\omega_k$.} \label{fig4a-1}
\end{figure}

\par By using the forward Euler method, the kinematic model of the AGV is reformulated as a discrete-time system \cite{kim2022lateral}, as shown in the following equations:
\begin{equation}
	\begin{cases}
		x_{k+1} = x_k + \Delta\cdot v_k{\rm cos}\theta_k , \\
		y_{k+1} = y_k + \Delta\cdot v_k{\rm sin}\theta_k,\\
		\theta_{k+1} = \theta_k + \Delta\cdot \omega_k,
	\end{cases}                         \label{4a-2}
\end{equation}
where $k=0,1,...,N-1$ represents the sampling times, and $\Delta$ represents the time step. 
The trajectory tracking task involves ensuring that the AGV follows a predetermined desired trajectory.
Define the following quadratic cost function
\begin{align}
	\nonumber J = &\sum_{k=0}^{N}\big[(X(k)-X_r(k))^TQ(X(k)-X_r(k)) 
	\\
	&\hspace*{2em}+ (U(k)-U_r(k))^TR(U(k)-U_r(k))\big],  \label{4a-3}
\end{align} 
where $X(k)=[x_k,y_k,\theta_k]^T$, $X_r(k)=[x_{r,k},y_{r,k},\theta_{r,k}]^T$, $U(k)$ $=[v_k,\omega_k]^T$, $U_r(k)$ $=[v_{r,k},\omega_{r,k}]^T$, $Q$ is a positive semi-definite matrix, and $R$ is a positive definite matrix. The task of trajectory tracking can be formulated as an optimal control problem, where the objective is to optimize linear velocity $v_k$ and angular $\omega_k$ to minimize the cost function (\ref{4a-3}) subject to the dynamical system (\ref{4a-2}) and the practical state or control constraints.

\subsection{Simulation Result}
The above trajectory tracking problem is addressed using MPC, as is shown in (\ref{3d-1})-(\ref{3d-3}). To evaluate the advantage of the proposed algorithm, the nonlinear optimization method are selected as benchmark. The algorithms are implemented using efficient optimization solvers, \texttt{fmincon} in MATLAB and CasADi, to solve the MPC problem \cite{worthmann2015model, li2021design}, employing SQP and interior point method.

\par The desired trajectory is generated by the system (\ref{4a-2}) with $v_{r,k}=2.3\text{m/s}$ and $\omega_{r,k}=0.0\text{rad/s}$.
The constraints are formulated as follows:
\begin{align}
	\nonumber &\hspace*{0.7em}\text{2.0} \text{m/s} \leq v_k \leq \text{2.35}\text{m/s},\ 
	\text{-1.5}\text{rad/s} \leq \omega_k \leq \text{1.0}\text{rad/s},
	\\
	\nonumber (x_k&-\text{3.0})^2+y_k^2 \geq \text{0.61}^2,
	(x_k-\text{6.1})^2+(y_k+\text{1.0})^2 \geq \text{0.81}^2,
	\\
	&\hspace*{3.5em}(x_k-\text{10.0})^2+(y_k-\text{0.4})^2 \geq \text{1.02}^2.   \label{4a-4}
\end{align}
The parameters in the above algorithms are selected as $\Delta$ $=\text{0.05}$, $[x_0,y_0,\theta_0]^T=[\text{0.0, -1.0, 0.0}]^T$, $N=\text{160}$, $N_p=\text{10}$, $Q=\text{diag}\{ \text{1.0, 1.0, 1.0} \}$, $R=\text{diag}\{ \text{1.1, 0.1} \}$, and $\mathcal{R}=\text{8.0}\boldsymbol{I}$ with $\boldsymbol{I}$ being the identity matrix of compatible dimensions.
The simulations are conducted on a computer equipped with a 13th Gen Intel(R) Core(TM) i7-13620H processor and 16.0 GB of RAM, running Matlab R2020b.
The tracking results and computational time comparisons for all algorithms are presented in Fig. \ref{fig4a-1} and Table \ref{table4a-1}, respectively.

\par Fig. \ref{fig4a-1} presents the state and control trajectories obtained by Algorithm 1 for solving the trajectory tracking problem (\ref{4a-2})-(\ref{4a-3}). In Fig. \ref{fig4a-1}(a), three circular obstacles, represented by black lines, are positioned at different locations. The green line depicts the tracking trajectory obtained using \texttt{fmincon} in the absence of obstacles, serving as a baseline for evaluating obstacle avoidance performance. The red line shows the tracking trajectory generated by Algorithm 1 in the presence of obstacles. The results demonstrate that Algorithm 1 effectively achieves obstacle avoidance.
Additionally, Fig. \ref{fig4a-1}(c) and (d) show the tracking performance of linear and angular velocities, respectively. The results indicate that the designed controller effectively follows the reference values while ensuring compliance with the specified constraints.

\par In Table 1, the overall and average computation times of Algorithm 1, \texttt{fmincon}, and CasADi are compared.
The total time is the average value obtained from 8 independent experiments for each algorithm.
Among these, the interior-point method in \texttt{fmincon} performs the least efficiently. 
The SQP in \texttt{fmincon} and the IPOPT solver in CasADi show similar performance, with average computation times under \text{10} ms. 
In contrast, Algorithm \text{1} is \text{4} to \text{5} times faster than the above two methods, further demonstrating its superior computational efficiency. 
The reason for this speed advantage can be found in Remark 3.

\section{Conclusion}
This paper addresses the discrete-time nonlinear optimal control problem under state and control constraints. 
By introducing an ALM to handle the constraints, a numerical solution algorithm is proposed.
The design of the optimization algorithm, along with the computation of gradient and Hessian matrix, revolves around analyzing and solving the corresponding FBDEs.
Notably, the entire process operates without reliance on external optimizers, enhancing computational efficiency and reducing deployment space requirements.

\bibliographystyle{IEEEtran}
\bibliography{reference_control_constraint}

\end{document}